\documentclass{amsart}

\usepackage{a4,xy,amssymb,diagrams}
\usepackage[hyperindex,pdfmark]{hyperref}

\xyoption{poly}
\xyoption{arc}
\xyoption{knot}
\xyoption{curve}
\xyoption{arrow}
\xyoption{all}

\makeindex

\newcommand{\dummy[1]}{{#1}}
\newcommand{\C}{\mathbb{C}}

\newcommand{\N}{\mathbb{N}}
\newcommand{\Z}{\mathbb{Z}}

\newcommand{\np}{\noindent}

\newcommand{\Oscr}{{\mathcal O}}

\newcommand{\Af}{\mathbb{A}}

\newcommand{\rootQ}{\sqrt[n]{\C Q}_{\dbrls \wisk{ab} \dbrrs}}

\newcommand{\dbrl}{[ \hskip -4pt [}
\newcommand{\dbrr}{] \hskip -4pt ]}
\newcommand{\dbrls}{\text{\tiny [ \hskip -2pt [}}
\newcommand{\dbrrs}{\text{\tiny ] \hskip -2pt ]}}

\newcommand{\kat}[1]{{\text{\bf \usefont{OT1}{pag}{m}{n} #1}}}

\newcommand{\wisk}[1]{{\text{\bf \usefont{OT1}{pag}{m}{n} #1}}}

\newcommand{\atn}{\text{\bf @}_{\text{\bf n}}}

\newtheorem{theorem}{Theorem}[section]

\newtheorem{proposition}[theorem]{Proposition}

\theoremstyle{definition}

\newtheorem{example}[theorem]{Example}

\theoremstyle{remark}

\title{Formal structures and representation spaces}

\author{Lieven Le Bruyn}
\address{Lieven Le Bruyn\\Universitaire Instelling Antwerpen\\B-2610 Antwerpen (Belgium)}
\email{lieven.lebruyn@ua.ac.be}
\author{Geert Van de Weyer}
\address{Geert Van de Weyer\\Universitaire Faculteiten Sint Ignatius Antwerpen\\B-2000 Antwerpen (Belgium)}
\email{geert.vandeweyer@ua.ac.be}

\begin{document}

\sloppy

\maketitle

\begin{abstract}
M. Kapranov introduced and studied in math.AG/9802041 the noncommutative formal structure of
a smooth affine variety. In this note we show that his construction is a particular case of
micro-localization and extend the construction functorially to representation schemes of
affine algebras. We describe explicitly the formal completions in the case of path algebras
of quivers and initiate the study of their finite dimensional representations.
\end{abstract}

\section{Introduction.}

Let $R$ be an associative
$\C$-algebra, $R^{Lie}$ its Lie structure and $R^{Lie}_m$ the subspace spanned by the
expressions $[r_1,[r_2,\hdots,[r_{m-1},r_m]\hdots ]$ containing $m-1$ instances of Lie brackets. The
{\it commutator filtration} of $R$ is the (increasing)
filtration by ideals $(F^k~R)_{k \in \Z}$ with $F^k~R = R$ for $d \in \N$
and
\[
F^{-k}~R = \underset{m}{\sum}~\underset{i_1 + \hdots + i_m = k}{\sum} R R^{Lie}_{i_1} R \hdots R R^{Lie}_{i_m} R 
\]
Observe that all $\C$-algebra morphisms preserve the commutator filtration.
The associated graded $gr_F~R$ is a (negatively) graded commutative Poisson algebra with part of degree
zero $R_{ab} = \tfrac{R}{[R,R]}$.

\np
Denote with $\kat{nil}_k$ the category of associative $\C$-algebras $R$ such that $F^{-k}R = 0$
(in particular, $\kat{nil}_1 = \kat{commalg}$ the category of commutative
$\C$-algebras). An algebra $A \in Ob(\kat{nil}_k)$ is said to be $k$-smooth if and only if for all
$T \in Ob(\kat{nil}_k)$, all nilpotent twosided ideals $I \triangleleft T$ and all $\C$-algebra
morphisms $A \rTo^{\phi} \tfrac{T}{I}$ there exist a lifted $\C$-algebra morphism
\[
\begin{diagram}
T & \rOnto & \dfrac{T}{I} \\
& \luDotsto_{\exists \tilde{\phi}} & \uTo^{\phi} \\
& & A
\end{diagram}
\]
making the diagram commutative. 

\np
Kapranov proves \cite[Thm 1.6.1]{Kapranov} that an affine commutative {\it smooth} algebra $C$ has a
{\it unique} (upto $\C$-algebra isomorphism identical on $C$) $k$-smooth {\it thickening} $C^{(k)}$ with 
$C_{ab}^{(k)} \simeq C$. The inverse limit (connecting morphisms are given by the uniqueness result)
\[
C^f = \underset{\leftarrow}{\text{lim}}~C^{(k)} \]
is then called the {\it formal completion} of $C$. Clearly, one has $C^f_{ab} = C$.

\begin{example}
Consider the affine space $\Af^d$ with coordinate ring
$\C[x_1,\hdots,x_d]$ and order the coordinate functions $x_1 < x_2 < \hdots < x_d$.
Let $\mathfrak{f}_d$ be the free Lie algebra on $\C x_1 \oplus \hdots \oplus \C x_d$ which has
an ordered basis $B = \cup_{k \geq 1} B_k$ defined as follows. $B_1$ is the ordered set
$\{ x_1,\hdots,x_d \}$ and $B_2 = \{ [x_i,x_j] \mid j < i \}$, ordered such that $B_1 < B_2$
and $[x_i,x_j] < [x_k,x_l]$ iff $j < l$ or $j=l$ and $i < k$. Having constructed the ordered
sets $B_l$ for $l < k$ we define 
\[
B_k = \{ [t,w] \mid t=[u,v] \in B_l, w \in B_{k-l} \text{\ such that \ } v \leq w < t \text{\ for \ }
l < k \}. \]
For $l < k$ we let $B_l < B_k$ and $B_k$ is ordered by $[t,w] < [t'.w']$ iff $w < w'$ or
$w=w'$ and $t < t'$. It is well known (see for example \cite[Ex. 5.6.10]{MKS}) that $B$ is an
ordered $\C$-basis of the Lie algebra $\mathfrak{f}_d$ and that its enveloping algebra
\[
U(\mathfrak{f}_d) = \C \langle x_1,\hdots,x_d \rangle \]
is the free associative algebra on the $x_i$. We number the elements of $\cup_{k \geq 2}B_k$
according to the order $\{ b_1,b_2,\hdots \}$ and for $b_i \in B_k$ we define $ord(b_i) = k-1$ (the
number of brackets needed to define $b_i$). Let $\Lambda$ be the set of all functions with finite
support $\lambda : \cup_{k \geq 2} B_k \rTo \N$ and define $ord(\lambda) = \sum \lambda(b_i) ord(b_i)$.
Rephrasing the Poincar\'e-Birkhoff-Witt result for $U(\mathfrak{f}_d)$ we have that any noncommutative
polynomial $p \in \C \langle x_1,\hdots,x_d \rangle$ can be written uniquely as a finite sum
\[
p = \sum_{\lambda \in \Lambda} \dbrl f_{\lambda} \dbrr~M_{\lambda} \]
where $\dbrl f_{\lambda} \dbrr \in \C[x_1,\hdots,x_d] = S(B_1)$ and $M_{\lambda} = \prod_i b_i^{\lambda(b_i)}$.
In particular, for every $\lambda,\mu,\nu \in \Lambda$, there is a unique bilinear differential operator
with polynomial coefficients
\[
C_{\lambda \mu}^{\nu} : \C[x_1,\hdots,x_d] \otimes_{\C} \C[x_1,\hdots,x_d] \rTo \C[x_1,\hdots,x_d] \]
defined by expressing the product $\dbrl f \dbrr~M_{\lambda}.~\dbrl g \dbrr~M_{\mu}$ in
$\C \langle x_1,\hdots,x_d \rangle$ uniquely as $\sum_{\nu \in \Lambda} \dbrl C_{\lambda \mu}^{\nu}(f,g)
\dbrr~M_{\nu}$, see \cite[Prop. 3.4.3]{Kapranov}. By associativity of $\C \langle x_1,\hdots,x_d \rangle$
the $C_{\lambda \mu}^{\nu}$ satisfy the associativity constraint, that is, we have equality of the
trilinear differential operators
\[
\sum_{\mu_1} C_{\mu_1 \lambda_3}^{\nu} \circ (C_{\lambda_1 \lambda_2}^{\mu_1} \otimes \wisk{id}) =
\sum_{\mu_2} C_{\lambda_1 \mu_2}^{\nu} \circ (\wisk{id} \otimes C_{\lambda_2 \lambda_3}^{\mu_2})
\]
for all $\lambda_1,\lambda_2,\lambda_3,\nu \in \Lambda$. Kapranov defines the algebra
$\C \langle x_1,\hdots,x_d \rangle_{\dbrls \wisk{ab} \dbrrs}$ to be the $\C$-vectorspace of possibly 
{\it infinite formal sums}
$\sum_{\lambda \in \Lambda} \dbrl f_{\lambda} \dbrr~M_{\lambda}$ with multiplication defined by the
operators $C_{\lambda \mu}^{\nu}$.

To prove that this algebra is the formal completion of $\C[x_1,\hdots,x_d]$ observe that the
commutator filtration
on $\C \langle x_1,\hdots,x_d \rangle$ has components
\[
F^{-k}~\C \langle x_1,\hdots,x_d \rangle = \{ \sum_{\lambda} \dbrl f_{\lambda} \dbrr~M_{\lambda} ,
\forall \lambda : ord(\lambda) \geq k \}
\]
Moreover,the quotient $\tfrac{\C \langle x_1,\hdots,x_d \rangle}{F^{-k}~\C \langle x_1,\hdots,x_d \rangle}$ is
$k$-smooth using the lifting property of free algebras and the fact that algebra morphisms preserve the
commutator filtration. Therefore,
\[
\C[x_1,\hdots,x_d]^f = \underset{\leftarrow}{\text{lim}}~\dfrac{\C \langle x_1,\hdots,x_d \rangle}{F^{-k}~\C \langle x_1,\hdots,x_d \rangle}
\simeq \C \langle x_1,\hdots,x_d \rangle_{\dbrls \wisk{ab} \dbrrs} . \]
\end{example}

\par \vskip 4mm
\noindent
Let $X$ be a (commutative) smooth affine variety (both assumptions are crucial !), then
Kapranov uses the formal completion of the algebra of functions, $\C[X]^f$, to define a sheaf of
noncommutative algebras on $X$, $\Oscr^f_X$ the Kapranov formal structure on $X$.

\np
From \cite[\S 4]{Kapranov} it follows that it is not possible in general to extend a Kapranov
formal structure on an arbitrary smooth variety $X$ from that on the affine open pieces. In fact, the
obstruction gives important new invariants of smooth varieties, related to Atiyah classes.

\np
When $X$ is affine, smoothness is essential to construct and prove uniqueness of the Kapranov formal
structure $\Oscr^f_X$. At present there is no natural functorial extension of formal structures
to arbitrary affine varieties. One of the major goals of this note is to construct such an
extension for representation spaces of affine non-commutative algebras.

\begin{example} Let us recall the construction of $\Oscr^f_{\Af^d}$ from \cite{Kapranov}.
Let $A_d(\C)$ be the $d$-th Weyl algebra, that is, the ring of differential operators with polynomial
coefficients on $\Af^d$. Let $\Oscr_{\Af^d}$ be the structure sheaf on $\Af^d$ then it is well-known that
the ring of sections $\Oscr_{\Af^d}(U)$ on any Zariski open subset $U \rInto \Af^d$ is a left
$A_d(\C)$-module. Kapranov defines in \cite[Thm. 3.5.3]{Kapranov} a sheaf $\Oscr^f_{\Af^d}$ of
noncommutative algebras on $\Af^d$ by taking as its sections over $U$ the
algebra
\[
\Oscr_{\Af^d}^f(U) = \C \langle x_1,\hdots,x_d \rangle_{\dbrls \wisk{ab} \dbrrs} 
\underset{\C[x_1,\hdots,x_d]}{\otimes} \Oscr_{\Af^d}(U) \]
that is the $\C$-vectorspace of possibly infinite formal sums $\sum_{\lambda \in \Lambda} \dbrl f_{\lambda} \dbrr~M_{\lambda}$ with
$f_{\lambda} \in \Oscr_{\Af^d}(U)$ and the multiplication is given as before by the action of the
bilinear differential operators $C_{\lambda \mu}^{\nu}$ on the left $A_d(\C)$-module $\Oscr_{\Af^d}(U)$,
that is, for all $f,g \in \Oscr_{\Af^d}(U)$ we have
\[
\dbrl f \dbrr~M_{\lambda}. \dbrl g \dbrr~M_{\mu} = \sum_{\nu} \dbrl C_{\lambda \mu}^{\nu}(f,g) \dbrr~M_{\nu}
\]
The sheaf of noncommutative algebras $\Oscr_{\Af^d}^f$ is called {\it Kapranov's formal structure} on
$\Af^d$.
\end{example}

\par
\vskip 4mm
\noindent
In the next section we will prove that the construction of Kapranov's formal structure
$\Oscr^f_X$ is a special case of micro-localization. In section three we will use this strategy to
extend the construction of formal structures in a functorial way to affine commutative schemes
$X = \wisk{rep}_n~A$, the scheme of all $n$-dimensional representations of an affine
$\C$-algebra $A$. In section four we will work out the special important case when $A$ is the path
algebra of a quiver and in the final section we will make some comments on the representation theory
of these formal completions.

\par \vskip 4mm

\noindent
{\bf Acknowledgement. } We like to thank F. Van Oystaeyen and A. Schofield for helpful discussions on the
contents of section 2, respectively 4.

\section{Microlocal interpretation.}

We recall briefly the algebraic construction of microlocalization. For more details we refer to the
monographs \cite{LiFVO} and \cite{FVObook}. Let $R$ be a filtered algebra with a separated
filtration $\{ F_n \}_{n \in \Z}$ and let $S$ be a multiplicatively closed subset of $R$ containing
$1$ but not $0$. For any $r \in F_n - F_{n-1}$ we denote its principal character
$\sigma(r)$ to be the image of $r$ in the associated graded algebra $gr(R)$. We assume that the
set $\sigma(S)$ is a multiplicatively closed subset of $gr(R)$ (always not containing $0$ whence
$\sigma$ is multiplicative on $S$). We define the {\it Rees ring} $\tilde{R}$ to be the graded
algebra
\[
\tilde{R} = \oplus_{n \in Z} F_n t^n \rInto R[t,t^{-1}] \]
where $t$ is an extra central variable. If $\sigma(s) \in gr(R)_n$ then we define the element
$\tilde{s} = st^n \in \tilde{R}_n$. The set $\tilde{S} = \{ \tilde{s}, s \in S \}$ is a multiplicatively
closed subset of homogeneous elements in $\tilde{R}$.

\np
Assume that $\sigma(S)$ is an Ore set in $gr(R) = \tfrac{\tilde{R}}{(t)}$, then for every $n \in \N_0$
the image $\pi_n(\tilde{S})$ is an Ore set in $\tfrac{\tilde{R}}{(t^n)}$ where $\tilde{R} \rOnto
\tfrac{\tilde{R}}{(t^n)}$ is the quotient morphism. Hence, we have an inverse system of graded
localizations and can form the inverse limit in the graded sense
\[
Q_{\tilde{S}}^{\mu}(\tilde{R}) = \underset{\leftarrow}{\text{lim}}^g~\pi_n(\tilde{S})^{-1} \dfrac{\tilde{R}}{(t^n)}
\]
The element $t$ acts torsionfree on this limit and hence we can form the filtered algebra 
\[
Q_S^{\mu}(R) = \dfrac{Q_{\tilde{S}}^{\mu}(\tilde{R})}{(t-1)Q_{\tilde{S}}^{\mu}(\tilde{R})} \]
which is the {\it micro-localization} of $R$ at the multiplicatively closed subset $S$. We recall that
the associated graded algebra of the microlocalization can be identified as
\[
gr(Q_S^{\mu}(R)) = \sigma(S)^{-1} gr(R). \]

Let $R$ be a $\C$-algebra with $R_{ab} = \tfrac{R}{[R,R]} = C$. We assume that the commutator filtration
$(F^k)_{k \in \Z}$ introduced in (1.3) is a separated filtration on $R$. Observe that this is not
always the case (for example consider $U(\mathfrak{g})$ for $\mathfrak{g}$ a semi-simple Lie
algebra) but often one can repeat the argument below replacing $R$ with $\tfrac{R}{\cap F^n}$.

\np
Observe that $gr(R)$ is a negatively graded commutative algebra with part of degree zero $C$. Take
a multiplicatively closed subset $S_c$ of $C$, then $S = S_c + [R,R]$ is a multiplicatively closed
subset of $R$ with the property that $\sigma(S) = S_c$ and clearly $S_c$ is an Ore set in $gr(R)$.
Therefore, $\tilde{S}$ is a multiplicatively closed set of the Rees ring $\tilde{R}$ consisting of
homogeneous elements of degree zero. Observing that $(t^n)_0 = F^{-n}t^n$ for all $n \in \N_0$ we see
that
\[
Q_S^{\mu}(R) = \underset{\leftarrow}{\text{lim}}~\pi_n(S)^{-1} \dfrac{R}{F^{-n}} \]
where $R \rOnto^{\pi_n} \tfrac{R}{F^{-n}}$ is the quotient morphism and $Q_S^{\mu}$ is filtered
again by the commutator filtration and has as associated graded algebra
\[
gr(Q^{\mu}_S(R)) = S_c^{-1} gr(R). \]
That is, the rings constructed in \cite[\S 2]{Kapranov} are just microlocalizations.

As explained in more detail in \cite{FVObook} one can define a {\it microstructure sheaf} $\Oscr_R^{\mu}$
on the affine scheme $X$ of $C$ by taking as its sections over the affine Zariski open set $X(f)$
\[
\Gamma(X(f),\Oscr_R^{\mu}) = Q^{\mu}_{S_f}(R) \]
where $S = \{ 1,f,f^2,\hdots \} + [R,R]$. Comparing with \cite{Kapranov} we have proved.

\begin{theorem} Let $X = \wisk{spec}~R_{ab}$ be smooth, then the microstructure sheaf $\Oscr_R^{\mu}$ coincides
with Kapranov's formal structure $\Oscr_X^f$.
\end{theorem}

An important remark to make is that one really needs microlocalization to construct a sheaf of
noncommutative algebras on $X$. If by some fluke we would have that all the $S_f$ are already
Ore sets in $R$, we might optimistically assume that taking as sections over $X(f)$ the Ore
localization $S_f^{-1} R$ we would define a sheaf $\Oscr_R$ over $X$. This is in general {\it not} the
case as the Ore set $S_g$ need no longer be Ore in a localization $S_f^{-1} R$ !

\section{Representation schemes.}

We will first show that representation schemes of noncommutative formally smooth algebras give
smooth affine varieties and are therefore endowed with a Kapranov formal structure.
A $\C$-algebra $A$ is said to be {\it formally smooth} if and only if it has the
lifting property with respect to {\it test-objects} in $\kat{alg}$. That is, let $I \triangleleft T$
be a nilpotent twosided ideal, then one can complete any $\C$-algebra morphism diagram
\[
\begin{diagram}
T & \rOnto & \dfrac{T}{I} \\
& \luDotsto_{\exists \tilde{\phi}} & \uTo_{\phi} \\
& & A
\end{diagram}
\]
If we restrict both $A$ and the testobjects $(T,I)$ to $\kat{commalg}$ we get Grothendieck's notion
of formally smooth commutative algebras. Grothendieck showed that if $A$ is essentially of finite type,
then this formal smoothness notion coincides with the usual geometric notion of smoothness.

\np
Motivated by this analogy, Quillen and Cunz \cite{CQ} have suggested to take formally smooth algebras as
coordinate rings of noncommutative smooth varieties.

\np
Let $A = \tfrac{\C \langle x_1,\hdots,x_d \rangle}{I_A}$ be an affine $\C$-algebra. The 
representation space $\wisk{rep}_n~A$
is the affine scheme representing the functor
\[
\begin{diagram}
\kat{commalg} & & & \rTo^{Hom_{\kat{alg}}(A,M_n(-))} & & \kat{sets} 
\end{diagram} . \]
Its coordinate ring $\C[\wisk{rep}_n~A]$ can be described as follows. For each $x_k$ consider an
$n \times n$ matrix of indeterminates $X_k = (x_{ij,k})_{i,j}$. For every relation
$f(x_1,\hdots,x_d) \in I_A$ consider the $n \times n$ matrix in $M_n(\C[x_{ij,k}, \forall i,j,k])
$ given by $f(X_1,\hdots,X_d) = (f_{ij})_{i,j}$. Each of the $f_{ij}$ is a polynomial in
$\C[x_{ij,k}, \forall i,j,k]$. Then we have
\[
\C[\wisk{rep}_n~A] = \dfrac{\C[x_{ij,k}, \forall i,j,k]}{(f_{ij}, \forall f \in I_A)} \]
Representability of the above functor equips us with a universal $\C$-algebra morphism
$A \rTo^{j_A} M_n(\C[\wisk{rep}_n~A])$ such that for every $\C$-algebra morphism $A \rTo^{\phi} M_n(C)$
with $C$ commutative
\[
\begin{diagram}
A & \rTo^{j_A} & & M_n(\C[\wisk{rep}_n~A]) \\
\dTo^{\phi} & \ldDotsto_{\exists ! M_n(\psi)} & & \\
M_n(C) & & &
\end{diagram}
\]
there is a unique morphism $\C[\wisk{rep}_n~A] \rTo^{\psi} C$,
making the above
diagram commute.

\np
Fore example, if $I_A=0$, that is $A \simeq \C \langle x_1,\hdots,x_d \rangle$, then
$\C[\wisk{rep}_n~A] = \C[x_{ij,k}, \forall i,j,k]$, whence
\[
\wisk{rep}_n~\C \langle x_1,\hdots,x_d \rangle = \Af^{dn^2} = \underbrace{M_n(\C) \oplus \hdots
\oplus M_n(\C)}_d . \]
The universal map $j_A$ maps $x_k$ to the {\it generic matrix} $X_k = (x_{ij,k})_{i,j}$ and the
morphism $\psi$ is determined by sending the variable $x_{ij,k}$ to the $(i,j)$-entry of the matrix 
$\phi(x_i)$.

\np
More generally, if $A$ is an affine formally smooth algebra, then all the representation spaces
$\wisk{rep}_n~A$ are affine smooth (commutative) varieties. Indeed, as $\C[\wisk{rep}_n~A]$ is affine,
it suffices to verify Grothendieck's formal smoothness property
\[
\begin{diagram}
T_c & \rOnto & \dfrac{T_c}{I} \\
& \luDotsto_{? \tilde{\phi}} & \uTo^{\phi} \\
& & \C[\wisk{rep}_n~A]
\end{diagram} \]
for $T_c$ a {\it commutative} algebra with nilpotent ideal $I$. Using formal smoothness of $A$ we have the
existence of $\psi$
\[
\begin{diagram}
M_n(T_c) & \rOnto & & M_n(\dfrac{T_c}{I}) \\
\uDotsto^{\exists \psi} & \luDotsto^{\exists M_n(\tilde{\phi})} & & \uTo_{M_n(\phi)} \\
A & \rTo_{j_A} & & M_n(\C[\wisk{rep}_n~A])
\end{diagram}
\]
and universality of $j_A$ provides us with the required lift $\tilde{\phi}$, proving 

\begin{proposition}
If $A$ is a formally smooth algebra then $\wisk{rep}_n~A$ is an affine smooth variety for all $n$.
\end{proposition}

For an arbitrary algebra $A$, however,
the representation space $\wisk{rep}_n~A$ is in general {\it not} smooth nor even reduced so as
mentioned before it is not immediately clear how to define a canonical and sufficiently
functorial formal structure on it. We will now show how this can be done.

The starting point is that for every associative algebra $A$ the functor
\[
\begin{diagram}
\kat{alg} & & \rTo^{Hom_{\kat{alg}}(A,M_n(-))} & &  \kat{sets} 
\end{diagram} \]
is {\it representable} in $\kat{alg}$. That is, there exists an associative $\C$-algebra
$\sqrt[n]{A}$ such that there is a natural equivalence between the functors
\[
Hom_{\kat{alg}}(A,M_n(-)) \underset{n.e.}{\sim} Hom_{\kat{alg}}(\sqrt[n]{A},-) . \]
In other words, for every associative $\C$-algebra $B$, there is a functorial one-to-one correspondence
between the sets
\[
\begin{cases}
\text{algebra maps} \quad A \rTo M_n(B)  \\
\text{algebra maps} \quad \sqrt[n]{A} \rTo B
\end{cases}
\]

\begin{example}
If $A = \C \langle x_1,\hdots,x_d \rangle$, then it is easy to see that
$\sqrt[n]{A} = \C \langle x_{11,1},\hdots,x_{nn,d} \rangle$. For, given an algebra map
$A \rTo^{\phi} M_n(B)$ we obtain an algebra map $\sqrt[n]{A} \rTo B$ by sending the free
variable $x_{ij,k}$ to the $(i,j)$-entry of the matrix $\phi(x_k) \in M_n(B)$. Conversely, to an
algebra map $\sqrt[n]{A} \rTo^{\psi} B$ we assign the algebra map $A \rTo M_n(B)$ by sending $x_k$
to the matrix $(\psi(x_{ij,k}))_{i,j} \in M_n(B)$. Clearly, these operations are each others inverses.
\end{example}

\par \vskip 4mm

\np
To define $\sqrt[n]{A}$ in general, consider the free algebra product $A \ast M_n(\C)$ and consider the subalgebra
\[
\sqrt[n]{A} = A \ast M_n(\C)^{M_n(\C)} = \{ p \in A \ast M_n(\C) \mid p.(1 \ast m) = 
(1 \ast m).p \ \forall m \in M_n(\C) \}
\]
Before we can prove the universal property of $\sqrt[n]{A}$ we need to recall a property that
$M_n(\C)$ shares with any Azumaya algebra : if $M_n(\C) \rTo^{\phi} R$ is an algebra morphism and
if $R^{M_n(\C)} = \{ r \in R \mid r.\phi(m) = \phi(m).r \ \forall m \in M_n(\C) \}$, then we have
$R \simeq M_n(\C) \otimes_{\C} R^{M_n(\C)}$. 

\np
In particular, if we apply this to $R = A \ast M_n(\C)$
and the canonical map $M_n(\C) \rTo^{\phi} A \ast M_n(\C)$ where $\phi(m) = 1 \ast m$ we obtain that
$M_n(\sqrt[n]{A}) = M_n(\C) \otimes_{\C} \sqrt[n]{A} = A \ast M_n(\C)$.

\np
Hence, if $\sqrt[n]{A} \rTo^{f} B$ is an algebra map we can consider the composition 
\[
A \rTo^{id_A \ast 1} A \ast M_n(\C) \simeq M_n(\sqrt[n]{A}) \rTo^{M_n(f)} M_n(B) \]
to obtain an algebra map $A \rTo M_n(B)$. Conversely, consider an algebra map $A \rTo^g M_n(B)$ and
the canonical map $M_n(\C) \rTo^i M_n(B)$ which centralizes $B$ in $M_n(B)$. Then, by the universal
property of free algebra products we have an algebra map $A \ast M_n(\C) \rTo^{g \ast i} M_n(B)$ and
restricting to $\sqrt[n]{A}$ we see that this maps factors
\[
\begin{diagram}
A \ast M_n(\C) & \rTo^{g \ast i} & M_n(B) \\
\uInto & & \uInto \\
\sqrt[n]{A} & \rDotsto & B
\end{diagram} \]
and one verifies that these two operations are each others inverses. 

\np 
It follows from the functoriality of the $\sqrt[n]{.}$ construction that $\C \langle x_1,\hdots,x_d \rangle \rOnto A$
implies that $\sqrt[n]{\C \langle x_1,\hdots,x_d \rangle} \rOnto \sqrt[n]{A}$. Therefore, if $A$ is
affine and generated by $\leq d$ elements, then $\sqrt[n]{A}$ is also affine and generated by
$\leq dn^2$ elements.

Next, we define a formal completion of $\C[\wisk{rep}_n~A]$ in a functorial way for any
associative algebra $A$.
Equip $\sqrt[n]{A}$ with the commutator filtration
\[
\begin{diagram}
\hdots & \rInto & F_{-2}~\sqrt[n]{A} & \rInto & F_{-1}~\sqrt[n]{A} & \rInto & \sqrt[n]{A} & = &
\sqrt[n]{A} & = & \hdots 
\end{diagram}
\]
Because algebra morphisms are commutator filtration preserving, it follows from the universal
property of $\sqrt[n]{A}$ that $\tfrac{\sqrt[n]{A}}{F_{-k} \sqrt[n]{A}}$ is the object in
$\kat{nil}_k$ representing the functor
\[
\begin{diagram}
\kat{nil}_k  & & \rTo^{Hom_{\kat{alg}}(A,M_n(-))}  & & \kat{sets} 
\end{diagram} . \]
In particular, because the categories $\kat{commalg}$ and $\kat{nil}_1$ are naturally 
equivalent, we deduce that
\[
\sqrt[n]{A}_{ab} = \dfrac{\sqrt[n]{A}}{[\sqrt[n]{A},\sqrt[n]{A} ]} = \dfrac{\sqrt[n]{A}}{F_{-1} \sqrt[n]{A}}
\simeq \C[\wisk{rep}_n~A] \]
because both algebras represent the same functor. We now define
\[
\sqrt[n]{A}_{\dbrls \wisk{ab} \dbrrs} = 
\underset{\leftarrow}{\text{lim}}~\dfrac{\sqrt[n]{A}}{F_{-k}~\sqrt[n]{A}} . \] 
Assume now that $A$ is formally smooth, then so is $\sqrt[n]{A}$ because we have seen before that
\[
M_n(\sqrt[n]{A}) \simeq A \ast M_n(\C) \]
and the class of formally smooth algebras is easily seen to be closed under free products and 
matrix algebras. Alternatively, one can apply Bergman's coproduct theorems \cite{Bergman} or
\cite[Thm. 2.20]{Schofield} for a strong version.

\np
As a consequence, we have for every $k \in \N$
that the quotient $\tfrac{\sqrt[n]{A}}{F_{-k} \sqrt[n]{A}}$ is $k$-smooth. Moreover, we have that
\[
(\dfrac{\sqrt[n]{A}}{F_{-k}~\sqrt[n]{A}})_{ab} \simeq \dfrac{\sqrt[n]{A}}{[\sqrt[n]{A},\sqrt[n]{A}]} 
 \simeq \C[\wisk{rep}_n~A] . \]
Because $\C[\wisk{rep}_n~A]$ is an affine commutative smooth algebra, we deduce from Kapranov's
uniqueness result of $k$-smooth thickenings that
\[
\C[\wisk{rep}_n~A]^{(k)} \simeq \dfrac{\sqrt[n]{A}}{F_{-k}~\sqrt[n]{A}} \]
and consequently that the formal completion of $\C[\wisk{rep}_n~A]$ can be identified with
\[
\C[\wisk{rep}_n~A]^f \simeq \sqrt[n]{A}_{\dbrls \wisk{ab} \dbrrs} . \]

\par \vskip 4mm

\begin{theorem}
Defining for an arbitrary $\C$-algebra $A$ the formal completion of $\C[\wisk{rep}_n~A]$ to be $\sqrt[n]{A}_{\dbrls \wisk{ab}
\dbrrs}$ gives a canonical extension of Kapranov's formal
structure on affine smooth commutative algebras to the class the coordinate rings of representation 
spaces on which it is functorial in the algebras. 
\end{theorem}

There is a natural action of $GL_n$ by algebra automorphisms on $\sqrt[n]{A}$. Let $u_A$ denote the
universal morphism $A \rTo^{u_A} M_n(\sqrt[n]{A})$ corresponding to the identity map on $\sqrt[n]{A}$.
For $g \in GL_n$ we can consider the composed algebra map
\[
\begin{diagram}
A & \rTo^{u_A} & M_n(\sqrt[n]{A}) \\
& \rdTo_{\psi_g} & \dTo_{g.g^{-1}} \\
& & M_n(\sqrt[n]{A})
\end{diagram}
\]
Then $g$ acts on $\sqrt[n]{A}$ via the automorphism $\sqrt[n]{A} \rTo^{\phi_g} \sqrt[n]{A}$ corresponding
the the composition $\psi_g$. It is easy to verify that this defines indeed a $GL_n$-action on
$\sqrt[n]{A}$.

\np
The formal structure sheaf $\Oscr^f_{\wisk{rep}_n~A}$ defined over $\wisk{rep}_n~A$, constructed from
$\sqrt[n]{A}$ by microlocalization as in the
foregoing section, will be denoted from now on by $\Oscr^f_{\sqrt[n]{A}}$. We see that it actually
has a $GL_n$-structure which is compatible with the $GL_n$-action on $\wisk{rep}_n~A$.

Finally we should clarify what representation theoretic information is contained in
$\sqrt[n]{A}_{\dbrls \wisk{ab} \dbrrs}$. The reduced variety of $\wisk{rep}_n~A$ gives information
about the $\C$-algebra maps $A \rTo M_n(\C)$. The scheme structure of $\wisk{rep}_n~A$ gives us the
$\C$-algebra maps $A \rTo M_n(C)$ where $C$ is a finite dimensional commutative $\C$-algebra. The
formal structure now gives us the $\C$-algebra maps $A \rTo M_n(B)$ where $B$ is a finite
dimensional noncommutative but basic $\C$-algebra. Recall that an algebra is said to be basic if all
its simple representations have dimension one.

\section{Path algebras of quivers.}

Even in the case when $A$ is formally smooth it is by no means easy to describe and
manipulate the $n$-th root algebra $\sqrt[n]{A}$ and its corresponding formal completion
$\sqrt[n]{A}_{\dbrls \wisk{ab} \dbrrs}$. In this and the next section we will discuss these
facts in the special (but important) case of path algebras of quivers.

Let $Q$ be a {\it quiver}, that is a directed graph on a finite set $Q_v = \{ v_1,\hdots,v_k \}$ of
vertices having a finite set $Q_a = \{ a_1,\hdots,a_l \}$ of arrows where we allow multiple
arrows between vertices and loops in vertices. We will depict vertex $v_i$ by \ \ 
$\xy 
\POS (0,0) *\cir<4pt>{}*+{\txt\tiny{i}} ="v2" \endxy$ \ \ and an arrow $a$ from vertex $v_i$ to
$v_j$ by \ \ $\xy 
\POS (0,0) *\cir<4pt>{}*+{\txt\tiny{j}} ="v2"
   , (15,0) *\cir<4pt>{}*+{\txt\tiny{i}} ="v3"
\POS"v3" \ar "v2"_{a} \endxy$.

\np
The path algebra $\C Q$ has as underlying $\C$-vectorspace basis the set of all oriented paths
in $Q$, including those of length zero corresponding to the vertices $v_i$. Multiplication in
$\C Q$ is induced by (left) concatenation of paths. More precisely, $1= v_1 + \hdots + v_k$ is a
decomposition of $1$ into mutually orthogonal idempotents and further we define
\begin{itemize}
\item{$v_j.a$ is always zero unless \ \  $\xy 
\POS (0,0) *\cir<4pt>{}*+{\txt\tiny{j}} ="v2"
   , (15,0) *\cir<4pt>{} ="v3"
\POS"v3" \ar "v2"_{a} \endxy$ \ \ in which case it is the path $a$,}
\item{$a.v_i$ is always zero unless \ \ $\xy 
\POS (15,0) *\cir<4pt>{}*+{\txt\tiny{i}} ="v2"
   , (0,0) *\cir<4pt>{} ="v3"
\POS"v2" \ar "v3"_{a} \endxy$ \ \ in which case it is the path $a$,}
\item{$a_i.a_j$ is always zero unless \ \ $\xy 
\POS (15,0) *\cir<4pt>{} ="v2"
   , (0,0) *\cir<4pt>{} ="v1"
   , (30,0) *\cir<4pt>{} ="v3"
\POS"v2" \ar "v1"_{a_i}
\POS"v3" \ar "v2"_{a_j} \endxy$ \ \ in which case it is the path $a_ia_j$.}
\end{itemize}
In order to see that $\C Q$ is formally smooth, take an algebra $T$ with a nilpotent twosided 
ideal $I \triangleleft T$ and consider
\[
\begin{diagram}
T & \rOnto & \dfrac{T}{I} \\
& \luDotsto_{? \tilde{\phi}} & \uTo_{\phi} \\
& & \C Q \end{diagram} \]
The decomposition $1 = \phi(v_1) + \hdots + \phi(v_k)$ into mutually orthogonal idempotents in
$\tfrac{T}{I}$ can be lifted up the nilpotent ideal $I$ to a decomposition $1 = \tilde{\phi}(v_1) +
\hdots + \tilde{\phi}(v_k)$ into mutually orthogonal idempotents in $T$. But then, taking for every
arrow $a$
\[
\xy 
\POS (0,0) *\cir<4pt>{}*+{\txt\tiny{j}} ="v2"
   , (15,0) *\cir<4pt>{}*+{\txt\tiny{i}} ="v3"
\POS"v3" \ar "v2"_{a} \endxy \quad \text{an arbitrary element} \quad \tilde{\phi}(a) \in
\tilde{\phi}(v_j) ( \phi(a) + I ) \tilde{\phi}(v_i) \]
gives a required lifted algebra morphism $\C Q \rTo^{\tilde{\phi}} T$.

Next, we will describe the smooth affine schemes $\wisk{rep}_n~\C Q$.
Consider the semisimple subalgebra
$V=\underbrace{\C \times \hdots \times \C}_k$ generated by the vertex-idempotents $\{ v_1,\hdots,v_k \}$.
Every $n$-dimensional representation of $V$ is semi-simple and determined by the multiplicities by
which the factors occur. That is, we have a decomposition
\[
\wisk{rep}_n~V = \bigsqcup_{\sum a_i = n} GL_n / (GL_{a_1} \times \hdots \times GL_{a_k}) =
\bigsqcup_{\alpha} \wisk{rep}_{\alpha}~V \]
into homogeneous spaces where $\alpha$ runs over the {\it dimension vectors} $\alpha = (a_1,\hdots,
a_k)$ such that $\sum_i a_i = n$. The inclusion $V \rInto \C Q$ induces a map $\wisk{rep}_n~\C Q \rTo^{\psi}
\wisk{rep}_n~V$ and we have the decomposition of $\wisk{rep}_n~\C Q$ into associated fiber bundles
\[
\psi^{-1}(\wisk{rep}_{\alpha}~V) = GL_n \times^{GL(\alpha)} \wisk{rep}_{\alpha}~Q . \]
Here, $GL(\alpha) = GL_{a_1} \times \hdots \times GL_{a_k}$ embedded along the diagonal in $GL_n$ and
$\wisk{rep}_{\alpha}~Q$ is the affine space of $\alpha$-dimensional representations of the quiver
$Q$. That is,
\[
\wisk{rep}_{\alpha}~Q = \bigoplus_{ \xy 
\POS (0,0) *\cir<4pt>{}*+{\txt\tiny{j}} ="v2"
   , (15,0) *\cir<4pt>{}*+{\txt\tiny{i}} ="v3"
\POS"v3" \ar "v2"_{a} \endxy } 
M_{a_j \times a_i}(\C) \]
and $GL(\alpha)$ acts on this space via base-change in the vertex-spaces. That is, $\wisk{rep}_n~\C Q$
is the disjoint union of smooth affine components depending on the dimension vectors $\alpha = (a_1,
\hdots,a_k)$ such that $\sum a_i = n$.

\par \vskip 4mm

The importance of the class of path algebras comes from the fact that if $A$ is an arbitrary affine
formally smooth algebra, the {\it \'etale local} $GL_n$-structure of the representation space
$\wisk{rep}_n~A$ is determined by a $\wisk{rep}_{\alpha}~Q$ for a certain quiver $Q$ and dimension
vector $\alpha$. We refer to \cite{LBngatn} or \cite{LBetale} for more
details.

\np
We have seen before that $\sqrt[n]{\C Q}$ is formally smooth, hence we would like to determine the
quiver settings relevant for the study of the representation space $\wisk{rep}_m~\sqrt[n]{\C Q}$
for arbitrary $m \in \N$. Before we can do this we need to recall the notion of {\it universal
localization}. We refer to \cite[Chp. 4]{Schofield} for full details.

\par
\vskip 4mm

\np
Let $A$ be a $\C$-algebra and $\kat{projmod}~A$ the category of finitely generated projective left $A$-modules.
Let $\Sigma$ be some class of maps in
this category (that is some left $A$-module morphisms between certain projective modules). In \cite[Chp. 4]{Schofield}
it is shown that there exists an algebra map $A \rTo^{j_{\Sigma}} A_{\Sigma}$ with the universal
property that the maps $A_{\Sigma} \otimes_A \sigma$ have an inverse for all $\sigma \in \Sigma$.
$A_{\Sigma}$ is called the universal localization of $A$ with respect to the set of maps $\Sigma$.

\np
When $A$ is formally smooth, so is $A_{\Sigma}$. Indeed, consider a test-object $(T,I)$ in
$\kat{alg}$, then we have the following diagram
\[
\begin{diagram}
T & \rOnto & \dfrac{T}{I} \\
\uDotsto_{\psi} & \luDotsto_{\tilde{\phi}} & \uTo^{\phi} \\
A & \rTo_{j_{\Sigma}} & A_{\Sigma}
\end{diagram}
\]
where $\psi$ exists by smoothness of $A$. By Nakayama's lemma all maps $\sigma \in \Sigma$ become
isomorphisms under tensoring with $\psi$. Then, $\tilde{\phi}$ exists by the universal property
of $A_{\Sigma}$.

\np
Consider the special case when $A$ is the path algebra $\C Q$ of a quiver on $k$ vertices. Then, we can
identify the isomorphism classes in $\kat{projmod}~\C Q$ with $\N^k$. To each vertex $v_i$ corresponds
an {\it indecomposable} projective left $\C Q$-ideal $P_i$ having as $\C$-vectorspace basis all paths
in $Q$ {\it starting at} $v_i$ (similarly, there is an indecomposable projective right $\C Q$-ideal
$P_i^r$ with basis the paths {\it ending} in $v_i$). We can also determine the space of homomorphisms
\[
Hom_{\C Q}(P_i,P_j) = \bigoplus_{\xy 
\POS (0,0) *\cir<4pt>{}*+{\txt\tiny{i}} ="v2"
   , (15,0) *\cir<4pt>{}*+{\txt\tiny{j}} ="v3"
\POS"v3" \ar@{~>} "v2"_{p} \endxy } \C p \]
where $p$ is an oriented path in $Q$ starting at $v_j$ and ending at $v_i$. Therefore, any
$A$-module morphism $\sigma$ between two projective left modules
\[
P_{i_1} \oplus \hdots \oplus P_{i_u} \rTo^{\sigma} P_{j_1} \oplus \hdots \oplus P_{j_v} \]
can be represented by an $u \times v$ matrix $M_{\sigma}$ whose $(p,q)$-entry $m_{pq}$ is a linear
combination of oriented paths in $Q$ starting at $v_{j_q}$ and ending at $v_{i_p}$.

\np
Now, form an $v \times u$ matrix $N_{\sigma}$ of free variables $y_{pq}$ and consider the algebra
$\C Q_{\sigma}$ which is the quotient of the free product $\C Q \ast \C \langle y_{11},\hdots,y_{uv}
\rangle$ modulo the ideal of relations determined by the matrix equations
\[
M_{\sigma}. N_{\sigma} = \begin{bmatrix} v_{i_1} & & 0 \\
& \ddots & \\
0 & & v_{i_u} \end{bmatrix} \quad \quad 
N_{\sigma} . M_{\sigma} = \begin{bmatrix}
v_{j_1} & & 0 \\ & \ddots & \\
0 & & v_{j_v} \end{bmatrix} \]
Equivalently, $\C Q_{\sigma}$ is the path algebra of a quiver {\it with relations} where the quiver
is $Q$ extended with arrows $y_{pq}$ from $v_{i_p}$ to $v_{j_q}$ for all $1 \leq p \leq u$ and
$1 \leq q \leq v$ and the relations are the above matrix entry relations.

\np
Repeating this procedure for every $\sigma \in \Sigma$ we obtain the universal localization
$\C Q_{\Sigma}$. Observe that if $\Sigma$ is a finite set of maps, then the universal localization
$\C Q_{\Sigma}$ is an affine algebra.

\par
\vskip 4mm

\np
It is easy to see that the representation space $\wisk{rep}_n~\C Q_{\sigma}$ is an affine Zariski
open subscheme (but possibly empty) of $\wisk{rep}_n~\C Q$. Indeed, if $m = (m_a)_a \in \wisk{rep}_{\alpha}~Q$,
then $m$ determines a point in $\wisk{rep}_n~\C Q_{\Sigma}$ if and only if the matrices
$M_{\sigma}(m)$ in which the arrows are all replaced by the matrices $m_a$ are invertible for all
$\sigma \in \Sigma$. In particular, this induces numerical conditions on the dimension
vectors $\alpha$ such that $\wisk{rep}_{\alpha}~Q_{\Sigma} \not= \emptyset$. Let $\alpha = (a_1,\hdots,a_k)$
be a dimension vector such that $\sum a_i = n$ then every $\sigma \in \Sigma$ say with
\[
P_1^{\oplus e_1} \oplus \hdots \oplus P_k^{\oplus e_k} \rTo^{\sigma} P_1^{\oplus f_1} \oplus
\hdots \oplus P_k^{\oplus f_k} \]
gives the numerical condition
\[
e_1 a_1 + \hdots + e_k a_k = f_1 a_1 + \hdots + f_k a_k . \]

\par
\vskip 4mm

Let $Q$ be a quiver on $k$ vertices and consider the extended quiver $\tilde{Q}_n$
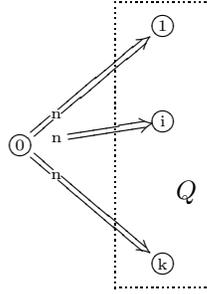
\[
\xy ;/r.15pc/:
(20,-30);
(20,30) **@{.}; (40,30) **@{.}; (40,-30) **@{.}; (20,-30) **@{.};
(35,-10) *+{Q};
\POS (30,25) *\cir<4pt>{}*+{\txt\tiny{1}} ="v1"
   , (30,5) *\cir<4pt>{}*+{\txt\tiny{i}} ="vi"
   , (30,-25) *\cir<4pt>{}*+{\txt\tiny{k}} ="vk"
   , (0,0) *\cir<4pt>{}*+{\txt\tiny{0}} ="v0"
   \POS"v0" \ar@{=>}|<>(.2)*{\txt{\tiny{\ n \ }}} "v1"
\POS"v0" \ar@{=>}|<>(.2)*{\txt{\tiny{\ n \ }}} "vi"
\POS"v0" \ar@{=>}|<>(.2)*{\txt{\tiny{\ n \ }}} "vk"
\endxy
\]
That is, we add to the vertices and arrows of $Q$ one extra vertex $v_0$ and for every vertex $v_i$ in $Q$ we add $n$ directed
arrows from $v_0$ to $v_i$. We will denote the $j$-th arrow $1 \leq j \leq n$ from $v_0$ to $v_i$ 
by $x_{ij}$.
From now on we will depict a bundle of $n$ arrows between two fixed vertices $v_k$ and $v_l$ by
$\xy
\POS (0,0) *\cir<4pt>{}*+{\txt\tiny{1}} ="v0"
   , (20,0) *\cir<4pt>{}*+{\txt\tiny{i}} ="vi"
\POS"v0" \ar@{=>}|<>(.2)*{\txt{\tiny{\ n \ }}} "vi"
\endxy$. Consider the morphism between projective left $\C \tilde{Q}_n$-modules
\[
P_1 \oplus P_2 \oplus \hdots \oplus P_k \rTo^{\sigma} \underbrace{P_0 \oplus \hdots \oplus P_0}_n \]
determined by the matrix
\[
M_{\sigma} = \begin{bmatrix} x_{11} & \hdots & \hdots & x_{1n} \\
\vdots & & & \vdots \\
x_{k1} & \hdots & \hdots & x_{kn} \end{bmatrix}. \]
We consider the universal localization $\C \tilde{Q}_{n \sigma}$, that is, we add for each vertex
$v_i$ in $Q$ another $n$ arrows $y_{ij}$ with $1 \leq j \leq n$ from $v_i$ to $v_0$
With these arrows $y_{ij}$ one forms the $n \times k$ matrix
\[
N_{\sigma} = \begin{bmatrix} y_{11} & \hdots & y_{k1} \\
\vdots & & \vdots \\
\vdots & & \vdots \\
y_{1n} & \hdots & y_{kn} \end{bmatrix} \]
and the universal localization $\C \tilde{Q}_{n \sigma}$ is described by the above
quiver with relations 
\[
M_{\sigma}.N_{\sigma} = \begin{bmatrix} v_1 & & 0 \\
& \ddots & \\ 0 & & v_k \end{bmatrix} \quad \text{and} \quad
N_{\sigma}.M_{\sigma} = \begin{bmatrix}
v_0 & & & 0 \\
& \ddots & & \\
& & \ddots & \\
0 & & & v_1 \end{bmatrix}. \]
We will depict this quiver with relations by the picture $\tilde{Q}_{n \sigma}$
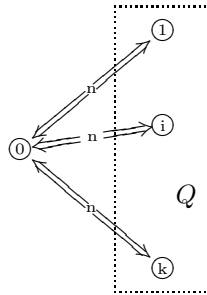
\[
\xy ;/r.15pc/:
(20,-30);
(20,30) **@{.}; (40,30) **@{.}; (40,-30) **@{.}; (20,-30) **@{.};
(35,-10) *+{Q};
\POS (30,25) *\cir<4pt>{}*+{\txt\tiny{1}} ="v1"
   , (30,5) *\cir<4pt>{}*+{\txt\tiny{i}} ="vi"
   , (30,-25) *\cir<4pt>{}*+{\txt\tiny{k}} ="vk"
   , (0,0) *\cir<4pt>{}*+{\txt\tiny{0}} ="v0"
\POS"v0" \ar@{<=>}|<>(.5)*{\txt{\tiny{\ n \ }}} "v1"
\POS"v0" \ar@{<=>}|<>(.5)*{\txt{\tiny{\ n \ }}} "vi"
\POS"v0" \ar@{<=>}|<>(.5)*{\txt{\tiny{\ n \ }}} "vk"
\endxy
\]
From the discussion above  we conclude :

\begin{theorem}
With notations as before, we have an isomorphism of $\C$-algebras
\[
\sqrt[n]{\C Q} \simeq v_0~\C \tilde{Q}_{n \sigma}~v_0 . \]
\end{theorem}

\begin{proof}
Indeed, the right hand side is generated by all the oriented cycles in $\tilde{Q}_{n \sigma}$
starting and ending at $v_0$ and is therefore generated by the $y_{ip}x_{iq}$ and the
$y_{ip}ax_{jq}$ where $a$ is an arrow in $Q$ starting in $v_j$ and ending in $v_i$. If we have an
algebra morphism
\[
\C Q \rTo^{\phi} M_n(B) \]
then we have an associated algebra morphism
\[
v_0~\C \tilde{Q}_{n \sigma}~v_0 \rTo^{\psi} B \]
defined by sending $y_{ip} a x_{jq}$ to the $(p,q)$-entry of the $n \times n$ matrix $\phi(a)$ and
$y_{ip}x_{iq}$ to the $(p,q)$-entry of $\phi(v_i)$. The defining relations among the $x_{ip}$ and
$y_{iq}$ introduced before imply that $\psi$ is indeed an algebra morphism.
\end{proof}

\begin{example}
Let us specialize to the case when $A = \C \langle a,b \rangle$, that is $A$ is the path algebra
of the quiver
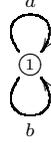
\[
\quad \xy 
\POS (0,0) *\cir<4pt>{}*+{\txt\tiny{1}} ="v1"
\POS"v1" \ar@(ul,ur)^a
\POS"v1" \ar@(dl,dr)_b
\endxy
\]
In order to describe $\sqrt[n]{A}$ we consider the quiver with relations
\[
\quad \xy 
\POS (0,0) *\cir<4pt>{}*+{\txt\tiny{1}} ="v1"
 , (-20,0) *\cir<4pt>{}*+{\txt\tiny{0}} ="v0"
\POS"v1" \ar@(ul,ur)^a
\POS"v1" \ar@(dl,dr)_b
\POS"v0" \ar@{<=>}|<>(.5)*{\txt{\tiny{\ n \ }}} "v1"
\endxy \quad  : \quad
y_i x_j = \delta_{ij} v_0 , \quad \quad \sum_i x_i y_i = v_1.
\]
We see that the algebra of oriented cycles in $v_0$ in this quiver with relations is
isomorphic to the free algebra in $2n^2$ free variables
\[
\C \langle y_1 a x_1, \hdots, y_n a x_n, y_1 b x_1, \hdots, y_n b x_n \rangle \]
which coincides with our knowledge of $\sqrt[n]{\C \langle a,b \rangle}$.
\end{example}

\section{Representations of $\sqrt[n]{\C Q}_{\dbrls \wisk{ab} \dbrrs}$.}

In this section we initiate the study of the finite dimensional representations of
$\rootQ$. The key observation is to observe that all simple representations of this
algebra are one-dimensional and that these one-dimensional representations are identified with
$\wisk{rep}_n~\C Q$ by definition of the root algebra. In order to describe all
$m$-dimensional representations of $\rootQ$ we look at the quotient map
\[
\wisk{rep}_m~\rootQ \rOnto^{\pi} \wisk{iss}_m~\rootQ \]
where $\wisk{iss}_m$ is the variety of isomorphism classes of $m$-dimensional semi-simple
representations (which is identified with the quotient variety of $\wisk{rep}_m$ under the
action of $GL_m$). First we will stratify this quotient variety by smooth subvarieties.

Any point $\xi$ of $\wisk{iss}_m~\rootQ$ corresponds to an $m$-dimensional semisimple representation of
the form
\[
M_{\xi} =S_1^{\oplus e_1} \oplus \hdots \oplus S_z^{\oplus e_z} \]
where the $S_i$ are distinct one-dimensional simple representations occurring with multiplicity $e_i$
and we choose the ordering of the components such that $e_1 \geq e_2 \geq \hdots \geq e_z > 0$.
That is, to $\xi$ we can associate a partition $\lambda(\xi) = (e_1,\hdots,e_z)$ of $m$. Every
partition $\lambda$ of $m$ determines a stratum $\wisk{iss}_m(\lambda)$ consisting of all points
$\xi$ such that $\lambda(\xi) = \lambda$. Using the
decomposition 
\[
\wisk{rep}_n~\C Q = \bigsqcup_{\alpha} GL_n \times^{GL(\alpha)} \wisk{rep}_{\alpha}~Q \]
(where $\alpha$ is a dimension vector of total dimension $n$) we conclude :

\begin{theorem} $\wisk{iss}_m~\rootQ$ has a stratification in locally closed smooth subvarieties
\[
\wisk{iss}_m(\lambda) \]
where $\lambda$ runs over the partition of $m$ and each stratum is the disjoint union of smooth
substrata
\[
\wisk{iss}_m(\lambda)(\alpha_1,\hdots,\alpha_z) \]
where $\alpha_i$ is a dimension vector for the quiver $Q$ of total dimension $n$ and $\lambda$ has $z$ nonzero parts.
\end{theorem}

Smoothness follows from the fact that if $\lambda = \lambda_1^{l_1} \hdots \lambda_y^{l_y}$, then
$\wisk{iss}_m(\lambda)$ is the product of $y$ smooth varieties, the $i$-th component being isomorphic
to the $l_i$-th symmetric power of $\wisk{rep}_n~Q$ with all the diagonals removed.

\par \vskip 4mm

Next we claim that the fibers of the quotient map $\pi$ are constant along these substrata
$\wisk{iss}_m(\lambda)(\alpha_1,\hdots,\alpha_k)$ can be identified with the nullcone of a specified quiver setting.
To prove this we relate $\wisk{rep}_n~\C Q$ to semistable representations of the extended quiver
$\tilde{Q}_n$.

If $S_i \in GL_n \times^{GL(\alpha_i)} \wisk{rep}_{\alpha_i}~Q$ (with $\alpha_i = (a_1,\hdots,a_k)$ we can use the $GL_n$-action to
assign to $S_i$ the representation $\tilde{S_i}$ of $\tilde{Q}_n$ of dimension vector
$\tilde{\alpha_i} = (1,\alpha_i)$ such that the restriction to the subquiver $Q$ is $S_i$, we put a
one-dimension space at vertex $v_0$ and take the arrows from $v_0$ to $v_j$ such that
\begin{itemize}
\item{the first $\sum_{l=1}^{j-1} a_l$ arrows map $\tilde{S}(v_0)$ to the zero vector.}
\item{the arrows $\sum_{l=1}^{j-1}a_l +1$ through $\sum_{l=1}^j a_l$ map $\tilde{S}(v_0)$ to a
standard basis of the vertexspace $S(v_j)$.}
\item{the remaining arrows all map $\tilde{S}(v_0)$ to the zero vector.}
\end{itemize}
If we take the character of $GL(\tilde{\alpha_i})$ determined by $\theta = (-m,1,\hdots,1)$ it follows
immediately that $\tilde{S_i}$ is a $\theta$-stable representation (that is, it contains no proper
subrepresentations of dimension vector $\tilde{\beta}$ such that $\theta(\tilde{\beta}) \leq 0$,
see \cite{King} for more details).

Hence, $M_{\xi}$ determines a $\theta$-semistable representation of $\tilde{Q}_n$. In
\cite{AL} it is proved that the \'etale local description of the corresponding moduli space
is determined by a local quiver setting. In our case the relevant quiver for 
the substratum determined by the partition $\lambda = (e_1,\hdots,e_z)$ and the dimension vectors
$\alpha_i$ is the  quiver setting $(\Gamma,\gamma)$ where $\Gamma$ is the quiver on $z$ vertices, say
$w_1,\hdots,w_z$ such that the number of directed arrows from $w_i$ to $w_j$ is equal to
\[
\delta_{ij} - \chi_{\tilde{Q}_n}(\tilde{\alpha_i},\tilde{\alpha_j}) \]
where $\chi_{\tilde{Q}_n}$ is the Euler form of the extended quiver, which is related to the
Euler form of the original quiver $Q$ by
\[
\chi_{\tilde{Q}_n} = \xy
\xymatrix@R-1.5pc@C-1.5pc{1 & -n & \dots & -n \\ 0 & \ar@{}[rrdd]|{\chi_Q}&&\\ \dots&&& \\ 0&&&
\save "2,2"."4,4"*+<-.1pc>[F--]\frm{}\restore}\drop\frm{(}\drop\frm{)}
\endxy \]
and the new dimension vector $\gamma = (e_1,\hdots,e_z)$ is determined by the multiplicities.

As we know that the Jordan-H\"older components of any representation in $\pi^{-1}(\xi)$ must be the
$S_i$, one deduces

\begin{theorem}
With notations as before, the fiber $\pi^{-1}(\xi)$ of the quotient morphism
$\wisk{rep}_m~\rootQ \rOnto^{\pi} \wisk{iss}_m~\rootQ$ in a point $\xi$ belonging to the substratum
$\wisk{iss}_m(\lambda)(\alpha_1,\hdots,\alpha_z)$ is determined by the nullcone of 
$\wisk{rep}_{\gamma}~\Gamma$. 
\end{theorem}

These observations are of particular importance in case the partition $\lambda = (m)$, that is if the
corresponding representation to $\xi$ is $M_{\xi} = S^{\oplus m}$. We then have an embedding
\[
\wisk{rep}_{\alpha}~Q \rInto^{i_m} M^{ss}_{\theta}(\tilde{Q}_n, m \tilde{\alpha}) \]
into the moduli space of $\theta$-semistable representations of $\tilde{Q}_n$ of dimension
vector $m.\tilde{\alpha}$.
The collection of these embeddings $i_m$ for $m \in \N$ is essentially equivalent to defining the
formal structure on $\wisk{rep}_{\alpha}~Q$. For more details we refer to \cite{LBncBS}.

\end{document}